\documentclass{elsart}
\usepackage{amssymb}
\usepackage{graphicx}

\begin{document}

\begin{frontmatter}

\title{A~bifurcational~geometric~approach to~the~problem~of~chaos~transition
in~the~classical~Lorenz~system\thanksref{label1}}
\thanks[label1]{This work was supported by the Netherlands Organization for Scientific
Research (NWO) and the German Academic Exchange Service (DAAD). The author is very grateful
to the Johann Bernoulli Institute for Mathematics and Computer Science of the University of
Groningen and the Institute for Mathematics of the TU Ilmenau for hospitality during his stays
in 2011--\,2012 and also to Henk Broer (Groningen) and Juergen Knobloch (Ilmenau) for very
fruitful discussions.
}

\author{Valery A. Gaiko}

\ead{valery.gaiko@gmail.com}

\address{National~Academy~of~Sciences~of~Belarus,~United~Institute~of~Informatics~Problems,
Surganov~Str.~6,~Minsk~220012,~Belarus}

\begin{abstract}
The classical Lorenz system is considered. For many years, this system has been the subject
of study by numerous authors. However, until now the structure of the Lorenz attractor is not
clear completely yet, and the most important question at present is to understand the bifurcation
scenario of chaos transition in this system. Using some numerical results and our bifurcational
geometric approach, we present a~new scenario of chaos transition in the classical Lorenz system.
\end{abstract}

\end{frontmatter}

\section{Introduction}
\label{1}

We consider a three-dimensional dynamical system
$$
\displaystyle\dot{x}=\sigma(y-x),\quad
\displaystyle\dot{y}=x(r-z)-y,\quad
\displaystyle\dot{z}= xy-bz
\eqno(1)
$$
known as the Lorenz system. Historically, (1) was the first dynamical system
for which the existence of an irregular attractor (chaos) was proved for
$\sigma=10,$ $b = 8/3,$ and $24,\!06<r<28.$ For many years, the Lorenz system
has been the subject of study by numerous authors; see, e.\,g.,
\cite{Sparrou1982}--\cite{Doedel2011}. However, until now the structure
of the Lorenz attractor is not clear completely yet, and the most important
question at present is to understand the bifurcation scenario of chaos transition
in system~(1).

In Section~2 of this paper, we recall the classical scenario of chaos transition
in the Lorenz system~(1). In Section~3, we give for (1) a relatively new chaos transition
scenario proposed by N.\,A.\,Magnitskii and S.\,V.\,Sidorov~\cite{MagnitskiiSidorov2006}.
In Section~4, we present a different bifurcation scenario for system~(1), where $\sigma=10,$
$b = 8/3,$ and $r>0,$ using numerical results of \cite{MagnitskiiSidorov2006} and our
bifurcational geometric approach to the global qualitative analysis of three-dimensional
dynamical systems which we applied earlier in the planar case \cite{Gaiko2003}--\cite{Gaiko12d}.

\section{The classical scenario of chaos tran\-si\-tion ($C$-scenario)}
\label{2}

First, let us briefly recall the contemporary point of view on the structure
of the Lorenz attractor and chaos transition \cite{Sparrou1982,MagnitskiiSidorov2006}.

1. The Lorenz system (1) is dissipative and symmetric with respect to the\linebreak $z$-axis.
The origin $O(0,0,0)$ is a singular point of system (1) for any $\sigma,$ $b,$ and~$r.$
It~is a stable node for $r<1.$ For $r=1,$ the origin becomes a triple singular point,
and then, for $r>1,$ there are two more singular points in the system:
$O_{1}(\sqrt{b(r-1)},\sqrt{b(r-1)},r-1)$ and $O_{2}(-\sqrt{b(r-1)},-\sqrt{b(r-1)},r-1)$
which are stable up to the parameter value $r_{a}=\sigma(\sigma+b+3)/(\sigma-b-1)$
$(r_{a}\approx24,\!74$ for $\sigma=10$ and $b=8/3).$ For all $r>1,$ the point $O$ is
a saddle-node. It~has a~two-dimensional stable manifold $W^{s}$ and a~one-dimensional
unstable manifold $W^{u}.$
If $1<r<r_{1}\approx13,\!9,$ then separatrices $\Gamma_{1}$ and $\Gamma_{2}$ issuing from
the point $O$ along its one-dimensional unstable manifold $W^{u}$ are attracted by their
nearest stable points $O_{1}$ and $O_{2},$ respectively.

2. If $r=r_{1},$ then each of the separatrices $\Gamma_{1}$ and $\Gamma_{2}$ becomes a closed
homoclinic loop. In this case, two homoclinic loops are tangent to each other and the $z$-axis
at the point $O$ and form a figure referred to as a homoclinic butterfly. It is assumed that
the generation of an unstable homoclinic butterfly is one of the two bifurcations leading to
the appearance of the Lorenz attractor.

3. If $r_{1}<r<r_{2}\approx24,\!06,$ then a saddle periodic trajectory bifurcates from each of
the closed homoclinic loops (these trajectories will be denoted by $L_{1}$ and $L_{2},$ respectively). In this case, separatrices $\Gamma_{1}$ and $\Gamma_{2}$ tend to the stable points
$O_{2}$ and $O_{1},$ respectively. It is usually assumed that stable manifolds of the saddle
periodic trajectories $L_{1}$ and $L_{2}$ are the boundaries of attraction domains of points
$O_{1}$ and $O_{2}.$ A curve issuing from the exterior of these domains can make oscillations
from the neighborhood of $L_{1}$ into a neighborhood of $L_{2}$ and conversely until it enters
the attraction domain of the attractor $O_{1}$ and $O_{2};$ the closer is parameter $r$
to the value $r_{2},$ the larger is the number of oscillations. This behavior of the system
is referred to as metastable chaos.
If $r=r_{2},$ then separatrices $\Gamma_{1}$ and $\Gamma_{2}$ do not tend to the points
$O_{2}$ and $O_{1},$ but wind around the limit saddle cycles $L_{2}$ and $L_{1},$ respectively.
Here the second bifurcation leading to the appearance of the Lorenz attractor takes place.
If $r_{2}<r<r_{3}=r_{a},$ then points $O_{1}$ and $O_{2}$ are still stable. In addition,
in the phase space, there is an attracting set $B$ referred to as the Lorenz attractor;
it is a set of integral curves moving from $L_{1}$ to $L_{2}$ and vice versa.
The saddle point $O,$ together with its separatrices $\Gamma_{1}$ and $\Gamma_{2},$
belongs to the attractor.

4. If $r\rightarrow r_{3}=r_{a},$ then the saddle limit cycles $L_{1}$ and $L_{2}$ shrink
to the points $O_{1}$ and $O_{2};$ for $r=r_{3},$ they vanish and coincide with these points
as a result of the Andronov--Hopf subcritical bifurcation.

5. If $r_{3}<r<r_{4}\approx30,\!1,$ then the Lorenz attractor is the unique stable limit set
of system (1). It is usually assumed that this set is a branching surface $S$ lying near
the plane $x-y=0$ and consisting of infinitely many sheets tied together and infinitely
close to each other. A phase trajectory issuing on the left from the $z$-axis comes untwisted
along a spiral around the point $O_{1}$ until the transition to the right of the $z$-axis,
after which it becomes untwisted along a spiral around the point $O_{2}$ in the opposite
direction. The number of rotations around the points $O_{1}$ and $O_{2}$ varies irregularly;
thus the motion looks chaotic. It is assumed that the attractor is not a two-dimensional
manifold and has a fractal structure \cite{MagnitskiiSidorov2006}.
If $r_{4}<r\lesssim313,$ then the structure of solutions of the system of Lorenz equations
becomes extremely complicated with alternation of chaotic and periodic modes. It is usually
assumed that there may be infinitely many periodicity windows in the system, and each of such
windows is a direct subharmonic cascade of bifurcations, which terminates with a basic stable
limit cycle. For further growth of $r,$ each of such cycles is destroyed by an intermittency,
and the appearance of periodicity windows is preceded by the inverse cascade of bifurcations
\cite{MagnitskiiSidorov2006}.

6. If $r>313,$ then the unique stable limit cycle is an attractor in the Lorenz system.

Thus, items 1--6 contain basic commonly accepted assertions dealing
with the Lorenz attractor and the scenarios of its appearance (and vanishing).
Note that all these assertions are based only on computer experiments
and speculative arguments rather than any analytic proofs. Some of these
assertions can readily be verified, and their validity is not brought into
question anywhere (e.\,g., the assertions of item~1).
However, other assertions are difficult to verify, and they always look quite
dubious. So, e.\,g., while saddle periodic cycles $L_{1}$ and $L_{2}$ are
really generated from homoclinic loops for $r=r_{1},$ and they are those determining
``eyes'' of the Lorenz attractor for $r=r_{2},$ but why are these ``eyes'' observed in
the attractor in the case $r>r_{3}$ as well in which the cycles $L_{1}$ and $L_{2}$
already vanished? The only possible conclusion is the following: the eyes of attractor
are not determined by the saddle cycles $L_{1}$ and $L_{2}$ even if they exist.
But if they also exist, at all it is unessential that they are born at $r=r_{1},$
as a result of bifurcation of a homoclinic butterfly. Also, assertions on the structure
of attractor and its dimension found on computer with incredible accuracy have always
been questionable. Finally, the phenomenon of intermittency did not find its logic
explanation. It was shown in \cite{MagnitskiiSidorov2006}, that actually in the Lorenz
system absolutely another scenario of chaos transition would be realized. We revise this
scenario too and, applying a bifurcational geometric approach, present a new scenario
of chaos transition in system (1) for $\sigma=10,$ $b = 8/3,$ and $r>0.$

\section{The Magnitskii--Sidorov scenario ($MS$-scenario)}
\label{3}

It turns out, see~\cite{MagnitskiiSidorov2006}, that all cycles from infinite family
of unstable cycles, generating Lorenz attractor, have crossing with an one-dimensional
unstable not invariant manifold $V^{u}$ of the point $O$ (do not confuse with the
invariant unstable manifold $W^{u}).$ This result follows from the theory of
dynamical chaos stated in \cite{MagnitskiiSidorov2006}. After the derivation
of analytic formulas for the manifold $V^{u},$ it becomes possible
to reduce the problem of establishing and proving the existence of
unstable cycles in the Lorenz system to the one-dimensional case, namely,
to finding stable points of the one-dimensional first return mapping defined
on the unstable manifold \cite{MagnitskiiSidorov2006}. By this method, it is
shown in \cite{MagnitskiiSidorov2006} that items~2 and~3 of the above-represented
classical scenario of transition to chaos in the Lorenz system~(1) are invalid.
Some assertions of items 4--6 fail, while other require a more detailed investigation.

1. This item remains the same as item~1 of the $C$-scenario.

2. If $r=r_{1}\approx13,\!9,$ then the separatrices $\Gamma_{1}$ and $\Gamma_{2}$
do not form two separate homoclinic loops. Here we have a bifurcation with the
generation of a single closed contour surrounding both stationary points
$O_{1}$ and $O_{2};$ the end of the separatrix $\Gamma_{1}$ enters the beginning
of the separatrix $\Gamma_{2},$ and vice versa, the end of $\Gamma_{2}$ enters
the beginning of $\Gamma_{1}.$ As $r$ grows, from this contour, a closed cycle
$C_{0}$ appears there first. It is an eight-shaped figure surrounding both points
$O_{1}$ and $O_{2}.$

3. If $r_{1}<r<r_{2}\approx24,\!06,$ then cycles $L_{1}$ and $L_{2}$ surrounding
the points $O_{1}$ and $O_{2},$ respectively, do not appear; but with further growth
of $r,$ pairs of cycles $C_{n}^{+},$ $C_{n}^{-},$ $n=0,1,\ldots,$ are successively
generated. They determine the generation of the Lorenz attractor. The cycle $C_{n}^{+}$
makes $n$ complete rotations in the half-space containing the point $O_{1}$ and one
incomplete rotation around the point $O_{2}.$ Conversely, the cycle $C_{n}^{-}$ makes
$n$~complete rotations around the point $O_{2}$ and one incomplete rotation around
the point $O_{1}.$

For each $r,$ $r_{1}<r<r_{2},$ there exists the number $n(r)$ $(n(r)\rightarrow\infty$
as $r\rightarrow~r_{2})$ such that in the phase-space of~(1), there are unstable cycles
$C_{0},$ $C_{k}^{+},$ $C_{k}^{-},$  $k=0,\ldots,n,$ and cycles $C_{km}^{+},$ $C_{km}^{-},$
$k,m<n,$ which make $k$ rotations around the point $O_{1}$ and $m$ rotations around the
point $O_{2}$ and are various combinations of the cycles $C_{n}^{+}$ and $C_{n}^{-},$ and
many other cycles generated by bifurcations of the cycles $C_{n}^{+}$ and $C_{n}^{-}$
\cite{MagnitskiiSidorov2006}.
Points of intersection of all these cycles with the manifold $V_{u}$ have the following
arrangement on the curve $V_{u}$ for $0\leq z_{min}\leq z\leq z_{max}<r-1$. The point
$z_{min}$ corresponds to the right large single loop of the cycle $C_{n}^{-}.$
This loop is the larger face of the right truncated cone of the set $S.$ Further,
the trajectory of the cycle passes into the left half-plane and makes $n$ clockwise
rotations around the point $O_{2}.$ The smallest first loop around the point $O_{2}$
is the smaller face of the truncated cone of the set $S.$ The point $z_{max}$
corresponds to the smallest loop of the cycle $C_{n}^{+}$ around the point $O_{1}.$
This loop is the smaller face of the right truncated cone. Further, the trajectory
of this cycle makes $n$ rotations around the point $O_{1}$ clockwise, passes
into the left half-plane, and makes one large rotation around the point $O_{2}.$
This rotation is the larger face of the left truncated cone. Between the points
$z_{min}$ and $z_{max}$ there is a point $z_{0}$ corresponding to the main
cycle $C_{0}.$

Boundaries of the attraction domains of the stable points $O_{1}$ and $O_{2}$
are given by the smallest loops of the cycles $C_{n}^{+}$ and $C_{n}^{-},$ whose
size decay as $r$ grows. Therefore, for some $r=r_{m},$ the attraction domain of
the set $B$ no longer intersects the attraction domains of points $O_{1}$ and
$O_{2},$ and the set $B$ becomes an attractor. Therefore, in the Lorenz system
$(a=10,$ $b=8/3),$ metastable chaos exists only in the interval $r_{1}<r<r_{m},$
and in the interval $r_{m}<r<r_{2},$ the system has three stable limit sets, namely,
$O_{1}$ and $O_{2}$ and the Lorenz attractor.

If $r\rightarrow r_{2},$ then the eye size decreases as the number of rotations of the
cycles $C_{n}^{+}$ and $C_{n}^{-}$ around the points $O_{1}$ and $O_{2},$ respectively,
grows. The value $z_{max}$ grows, and $z_{min}$ decays; moreover, $z_{min}\rightarrow0$
as $r\rightarrow r_{2}.$ The lengths of generatrices of truncated cones grow, since
additional rotations are added to the cone vertex and diminish the size of the
smaller face. Conversely, the larger face grows. If $r=r_{2},$ then $z_{min}=0,$
but $z_{max}<r-1;$ thus, the larger face of each cone achieves its maximal size,
while the smaller face is not contracted into a point, the cone vertex. The following
bifurcation takes place. In the limit as $n\rightarrow\infty,$ each set of cycles
$C_{n}^{+}$ (respectively, $C_{n}^{-})$ forms a point-cycle heteroclinic structure
consisting of two separatrix contours of the point $O.$ The first contour consists of
a separatrix issuing from the point $O$ along its unstable manifold and spinning on
the appearing (only for\linebreak $r=r_{2})$ saddle cycle $L_{1}$ (respectively, $L_{2})$
of the point $O_{1}$ (respectively, $O_{1})$. The second contour consists of the
separatrix spinning out from the saddle cycle $L_{1}$ (respectively,~$L_{2})$
and entering the point $O$ along its stable manifold.

As mentioned above, the described bifurcation does not lead to generation of
the Lorenz attractor for $r=r_{2}.$ It is more correct to say that it is only
a prerequisite of destruction of the attractor as $r$ decays. The attractor itself,
existing in the system for $r=r_{2},$ is formed from finitely many stable cycles
$C_{k}^{\pm},$ $k=0,\ldots,l,$ for $r<313.$ It contains neither separatrices
$\Gamma_{1}$ and $\Gamma_{2}$ of the point $O$ nor infinitely many unstable cycles
$C_{n}^{\pm}$ existing in the neighborhood of the point-cycle heteroclinic structure.

If $r_{2}<r<r_{3},$ then points $O_{1}$ and $O_{2}$ are still stable, and their
attraction domains are bound by the appearing limit cycles $L_{1}$ and $L_{2}$
contracting to points as $r\rightarrow r_{3}.$ But the Lorenz attractor $B$ is not
a set of integral curves going from $L_{1}$ to $L_{2}$ and back, and separatrices
$\Gamma_{1}$ and $\Gamma_{2}$ of the saddle point $O$ do not belong to the attractor.
Cycles $L_{1}$ and $L_{2}$ have already made their job at $r=r_{2}$ and no longer have
anything to do with the attractor. If $r_{2}<r<r_{3},$ then, just as in the case of
$r_{1}<r<r_{2},$ the cycles $C_{n}^{+}$ and $C_{n}^{-}$ appear again from separatrix
contours. The attractor is determined by finitely many such cycles \cite{MagnitskiiSidorov2006}.

4. For $r=r_{3},$ the saddle cycles $L_{1}$ and $L_{2}$ disappear. In the system,
there is a unique limit set, namely, the Lorenz attractor.

5. There exist one more important value of the parameter $r$ which affects
the formation of the Lorenz attractor. This is a point $r_{4}\approx30,\!485.$
If $r$ grows from $r_{3}$ to $r_{4},$ then the number of rotations of the cycles
$C_{n}^{+}$ and $C_{n}^{-}$ first rapidly decays, then grows again. In this case,
eyes by separatrices of the point $O$ are much smaller than attractor eyes and begin
to grow as $r$ increases. Therefore, $r_{4}$ is the point of minimum distance from
the line $(a=10,$ $b=8/3)$ in the space of parameters $(a,b,r)$ to the curve of
heteroclinic contours joining the point $O$ with the points $O_{1}$ and $O_{2}.$
The separatrices of the point $O$ approach one-dimensional stable manifolds of the
points $O_{1}$ and $O_{2}$ by the minimal distance but do not hit these points.
Therefore, almost heteroclinic and almost homoclinic contours exist in system~(1)
at the point~$r_{4}.$

The process of generation of the Lorenz attractor in system~(1) as $r$ decays from
the value $313$ up to $r_{4}$ is referred to as the incomplete double homoclinic
cascade \cite{MagnitskiiSidorov2006}. The complete cascade occurs if the $r$-axis
passes exactly through the point of existence of two homoclinic contours. Note that
in systems with a single homoclinic contour, there can be a simple complete or incomplete
homoclinic cascade of bifurcations of transition to chaos, and in \cite{MagnitskiiSidorov2006},
a detailed description of transition to chaos through the double homoclinic (complete or
incomplete) cascade of bifurcations is given. Just as in item~6 of the classical scenario,
if $r>313,$ then in the system, there exists a unique stable limit cycle $C_{0}$ surrounding
both points. If $r\approx313,$ then the cycle $C_{0}$ becomes unstable and generates
two stable cycles $C_{0}^{+}$ and $C_{0}^{-}$ which also surround the points $O_{1}$
and $O_{2}$ but have deflections in the direction of corresponding halves of
the unstable manifold $V^{u}$ of the point $O.$ This is the point where
the double homoclinic cascade of bifurcations really begins. In case of
an incomplete cascade, it consists of finitely many stages of appearance of stable cycles
$C_{k}^{\pm},$ $k=0,\ldots,l,$ and their infinitely many further bifurcations. But in case
of a complete cascade, the number of stages is infinite, and at the limit of $l\rightarrow\infty,$ cycles tend to homoclinic contours of the points $O_{1}$ and $O_{2},$ respectively. At~the $k$-th
stage of the cascade, originally stable cycles $C_{k}^{\pm}$ undergo a subharmonic cascade of
bifurcations and form two band-form attractors that consist of infinitely many unstable limit
cycles intersecting the respective domains of the unstable manifold $V^{u}$ of the point $O.$
Then these two bands merge and form a single attractor surrounding both the points $O_{1}$ and
$O_{2},$ after which there is a cascade of bifurcations of cycles generated as a result of the
merger and making rotations separately around the points $O_{1}$ and $O_{2}$ and simultaneously
around both the points. The last cascade of bifurcations has the property of self-organization,
since it is characterized by simplification of the structure of cycles and the generation of
new stable cycles with a smaller number of rotations around the points $O_{1}$ and $O_{2}$
as $r$ decays. Each cycle of the cascade of self-organization bifurcations undergoes its own
subharmonic cascade of bifurcations, after which all cycles formed during infinitely many
bifurcations of all subharmonic cascades and cascades of self-organization bifurcations of
cycles become unstable and form some set $B_{k}.$ After an incomplete homoclinic cascade of
bifurcations, we obtain a set $B=\bigcup\,B_{k}$ consisting of infinitely many possible unstable
cycles appearing at all stages of the cascade. These cycles generate an incomplete double
homoclinic attractor, that is the classical Lorenz attractor.

6. This item remains the same as item~6 of the $C$-scenario.

\section{The bifurcational geometric scenario ($G$-scenario)}
\label{4}

Revising the above scenarios, we present a new scenario of chaos transition in the Lorenz
system~(1) for $\sigma=10,$ $b = 8/3,$ and $r>0.$

1. If $r<1,$ the unique singular point $O$ of system~(1) is a stable node. For $r=1,$
it becomes a triple singular point, and then, for $r>1,$ there are two more singular
points in the system: $O_{1}$ and $O_{2}$ which are stable up to the parameter value
$r_{a}\approx24,\!74.$ For all $r>1,$ the point~$O$ is a saddle-node. It~has
a~two-dimensional stable manifold $W^{s}$ and an~one-dimensional unstable manifold $W^{u}.$
If $1<r<r_{l}=r_{1}\approx13,\!9,$ then the separatrices $\Gamma_{1}$ and $\Gamma_{2}$
issuing from the point $O$ along its one-dimensional unstable manifold $W^{u}$ are
attracted by their nearest stable points $O_{1}$ and $O_{2},$ respectively.

2. If $r=r_{l},$ then each of the separatrices $\Gamma_{1}$ and $\Gamma_{2}$ becomes
a closed homoclinic loop. In this case, two unstable homoclinic loops, $C_{0}^{+}$ and
$C_{0}^{-},$ are formed around the points $O_{1}$ and $O_{2},$ respectively. They are
tangent to each other and the $z$-axis at the point $O$ and form together a homoclinic
butterfly.

3. If $r_{l}<r<r_{a}\approx24,\!74,$ then, unfortunately, neither the $C$-scenario
nor the $MS$-scenario can be realized. The reason is that, in both cases, trajectories of
system~(1) should intersect the two-dimensional stable manifold $W^{s}$ of the point $O.$
Since this is impossible, the only way to overcome the contradiction is to suppose that
a~cascade of period-doubling bifurcations \cite{MagnitskiiSidorov2006} will begin immediately
in each of the half-spaces with respect to the manifold $W^{s},$ when $r>r_{l}.$ In this case,
each of the homoclinic loops $C_{0}^{+}$ and $C_{0}^{-}$ generates an unstable
limit cycle of period~2 which makes one rotation around the point $O_{1}$
and one rotation around the point $O_{2}$ but in the corresponding half-spaces
containing the points $O_{1}$ and $O_{2},$ respectively, and a stable limit cycle
of period~1 lying between the coils of the cycle of period~2. With further growth of
$r,$ each of the cycles of period~2 generates an unstable limit cycle of period~4 with
a stable limit cycle of period~3 inside of it and each of the cycles of period~1 generates
a stable limit cycle of period~2 with an unstable limit cycle of period~1 inside of it.
Then, after next doubling, we will have in each of the half-spaces an unstable limit cycle of
period~8 with an inserted stable limit cycle of period~7 and a stable limit cycle of period~6
with an inserted unstable limit cycle of period~5, and a stable limit cycle of period~4 with
an inserted unstable limit cycle of period~3, and an unstable limit cycle of period~2 with
an inserted stable limit cycle of period~1. Continuing this process further, we will obtain
limit cycles of all periods from one to infinity, and the space between these cycles will be
filled by spirals issuing from unstable limit cycles and tending to stable limit cycles as
$t\rightarrow+\infty.$ These cycles are inserted into each other, they make various combinations
of rotation around the points $O_{1}$ and $O_{2}$ in the corresponding half-spaces containing
these points and form geometric constructions (limit periodic sets) which look globally like
very flat truncated cones described in item~3 of the $MS$-scenario \cite{MagnitskiiSidorov2006}.

4. For $r=r_{a}\approx24,\!74,$ the biggest unstable limit cycles of infinite period disappear
through the Andronov--Shilnikiv bifurcation \cite{Shilnikov2001,Kuznetsov2004} in each of the
half-spaces containing the points $O_{1}$ and $O_{2}$ (the cone vertices are at these points),
and these points become unstable saddle-foci.

5. If $r_{a}<r<+\infty,$ then a~cascade of period-halving bifurcations
\cite{MagnitskiiSidorov2006} occurs in each of the half-spaces with respect to the manifold
$W^{s}.$ We have got again two symmetric with respect to the $z$-axis limit periodic sets
consisting of limit cycles of all periods which are inserted into each other and make
various combinations of rotation around the points $O_{1}$ and $O_{2}$ in the corresponding
half-spaces containing these points, and the space between the cycles is filled by spirals
issuing from unstable limit cycles and tending to stable limit cycles as $t\rightarrow+\infty.$
The biggest limit cycles of these sets are stable now, and with further growth of $r,$
the period-halving process makes them and the whole limit periodic sets more and more flat.
The obtained geometric constructions are the only stable limit sets of system~(1).
The spirals of the unstable saddle-foci $O_{1}$ and $O_{2}$ and the trajectories issuing from
infinity tend to these limit periodic sets (more precisely, to their stable limit cycles)
as $t\rightarrow+\infty.$ Just these stable limit periodic sets form two symmetric parts
of the so-called Lorenz attractor, and this really looks very chaotic.

6. If $r\rightarrow+\infty$ (numerically, when $r\gtrsim313),$ then the period-halving process
will be finishing and system~(1) will have two stable limit cycles in two phase half-spaces
containing the unstable saddle-foci $O_{1}$ and $O_{2}$ of~(1). This completes our scenario
of chaos transition in the Lorenz system~(1).

\end{document}